\documentclass[runningheads]{llncs}
\usepackage{amsmath}
\usepackage{amssymb}
\usepackage[applemac]{inputenc}
\usepackage[colorlinks]{hyperref}
\usepackage{url}
\usepackage{graphicx,epstopdf}
\usepackage{subcaption}
\usepackage[T1]{fontenc}
\usepackage{graphicx}
\begin{document}
\title{Temporal Alignment of Human Motion Data: \\ A Geometric Point of View \thanks{Supported by FWF grant I 5015-N, Institut CNRS Pauli, by grant F77 of the Austrian Science Fund FWF (SFB "Advanced Computational Design", SP5), TU Wien and University of Lille}}
%
%
\author{Alice~Barbora~Tumpach\inst{1,2}\orcidID{0000-0002-7771-6758} \and
Peter~K\'an \inst{3}\orcidID{0000-0001-7437-9955} }
%
\authorrunning{Tumpach and K\'an}
%
\institute{Institut CNRS Pauli, Oskar-Morgenstern-Platz 1, 1090 Vienna, Austria \and
University of Lille, Cit\'e scientifique, 59650 Villeneuve d'Ascq, France
\email{alice-barbora.tumpach@univ-lille.fr}\\
\url{http://math.univ-lille1.fr/~tumpach/Site/home.html} \and
Institute of Visual Computing and Human-Centered Technology, TU Wien, Vienna, Austria \\
\email{peterkan@peterkan.com}}

\maketitle              
\begin{abstract}
Temporal alignment is an inherent task in most applications dealing with videos: action recognition, motion transfer, virtual trainers, rehabilitation, etc. In this paper we dive into the understanding of this task from a geometric point of view: 
in particular, we show that the basic properties that are expected from  a temporal alignment procedure imply that the set of aligned motions to a template form a slice to a principal fiber bundle for the group of temporal reparameterizations. A temporal alignment procedure provides a reparameterization invariant projection onto this particular slice. This geometric presentation allows to elaborate a consistency check for testing the accuracy of any temporal alignment procedure. We give examples of alignment procedures from the literature  applied to motions of tennis players. Most of them use dynamic programming to compute the best correspondence between two motions relative to a given cost function. This step is computationally expensive (of complexity $O(NM)$ where $N$ and $M$ are the numbers of frames). Moreover most methods use features that are invariant by translations and rotations in $\mathbb{R}^3$, whereas most actions are only invariant by translation along and rotation around the vertical axis, where the vertical axis is aligned with the gravitational field. The discarded information contained in the vertical direction is crucial for accurate synchronization of motions. We propose to  incorporate keyframe correspondences into the dynamic programming algorithm based on coarse information extracted from the vertical variations, in  our case from the elevation of the arm holding the racket. The temporal alignment procedures produced are not only more accurate, but also computationally more efficient.

\end{abstract}

\section{Introduction}

This work deals with temporal alignment of motions and therefore is connected to motion analysis in general. 
A state-of-the-art report on motion similarity modeling was presented by Sebernegg et al. \cite{Sebernegg} and by Senin~\cite{Senin}.  We will use a geometric formulation based on group actions and invariances, that can also be used for registration tasks in Shape analysis. In the present setting, the relevant features are curves in $\mathbb{R}^3$ or in some Lie groups and homogeneous spaces. Curves in homogeneous spaces where used for  animation purposes by Celledoni et al. 
\cite{Celledoni}.  We improve the implementations of temporal alignment procedures introduced in previous research  \cite{Srivastava,Kacem,Celozzi,Celledoni2,Celledoni3}.
Here our main concern is to align motions in order to be able to display them in a synchronized manner. 

\subsection*{Contributions}
\begin{itemize}
\item[$\bullet$] we give a geometric formulation of the task consisting of temporal alignment of two motions  (Section~\ref{geometric});
\item[$\bullet$] this mathematical formulation provides a guide to any temporal alignment procedure, in particular we explain how the consistency of a temporal alignment procedure can be checked (Section~\ref{Consistency_check_section});
\item[$\bullet$] we compare different alignment procedures including a coarse alignment procedure that provides a simple and computationally efficient solution to which any other method should be compared in the effort of finding a balance between accuracy and complexity (see Geometric Green Learning~\cite{GGL} and Section~\ref{procedures});
\item[$\bullet$] we provide a variant of Dynamic Time Wraping algorithm that takes into account keyframes correspondences, is computationally more efficient and could be used for different purposes (Section~\ref{incorporating}).
\end{itemize}

\section{Temporal Alignment from a geometric point of view}\label{geometric}

\subsection{Type of data under consideration}

In this paper, we give a geometrical picture of the task consisting of temporal alignment of two motions. As a running example we describe the temporal alignment of motion data of tennis players performing the same action. We use the time evolution of extracted skeletons to characterize the motions. Skeletons extraction is performed using different devices (kinect, simi), but in our examples, the number of joints is fixed, and the extracted features are the same. 
A priori, the time interval in which a motion is performed could depend on the motion under consideration, but  as a pre-processing step we renormalize all motions to the time interval $[0, 1]$, i.e. any action start at  time $t= 0$ and ends at time $t = 1$. Nevertheless, the number of frames may differ from actions to actions, depending on the devices used to record the movement. This means that, for a given action,  we have a discrete set of times $(t_1, \dots, t_{K}) $,  where $K$ is the number of frames, at which we know the exact positions of  joints. 
\subsection{What is Temporal Alignment?}

The task of aligning two motions $M_1, M_2: [0, 1]\rightarrow \mathbb{R}^{3N}$ can be understood at the theoretical level as the task of finding an optimal time warping (in other words diffeomorphism) $\varphi: [0, 1] \rightarrow [0,1]$ such that $M_2\circ \varphi$ is visually as close as possible to $M_1$. On a practical level, the output of time alignment of motions $M_1$ and $M_2$ will be a correspondence between the set of frames of $M_1$, labelled by $\{1, \dots, K_1\}$ with the set of frames of $M_2$, labelled by $\{1, \dots, K_2\}$, where a given frame  can be in correspondence with multiple frames. Given an optimal time warping $\varphi: [0, 1] \rightarrow [0,1]$ between two motions, the computation of the correspondence between frames  presents no difficulty. Therefore, the main challenging task in temporal alignment of motions resides in the computation of an optimal diffeomorphism aligning one motion to the other.

\subsection{Temporal Alignment procedures as maps from motions to the group of diffeomorphisms}
To summarize the setting mathematically at this stage, we consider the manifold $\mathcal{C}$ of all smooth curves parameterized by $[0, 1]$ and with values in $\mathbb{R}^{3N}$:
\begin{equation}
\mathcal{C}  = \{ c: [0, 1]\rightarrow \mathbb{R}^{3N}, c \textrm{  smooth}\},
\end{equation}
 and we define the set $\mathcal{M}$ of motions of a skeleton with $N$ joints, fixed set of links $\mathcal{L} \subset \{1, \dots, N\}\times\{1, \dots, N\}$ and fixed lengths of bones $c_{jk}>0$, $\{j, k\}\in \mathcal{L}$  by
\begin{equation}
\mathcal{M} := \{ f \in\mathcal{C},
\| f_j(t) - f_k(t) \| = c_{jk}, \forall t \in [0, 1], \forall \{j, k\}\in\mathcal{L}\}.
\end{equation}

The set of time warpings is the group $\operatorname{Diff}^{+}([0, 1])$ of (orientation-preserving) diffeomorphisms $\varphi: [0, 1] \rightarrow [0, 1]$, (i.e. sending $0$ to $0$ and $1$ to $1$) acting on curves in $\mathcal{C}$, hence also on motions in $\mathcal{M}$, by
$
\varphi\cdot f := f\circ \varphi.
$
Given a reference motion $M_{\textrm{ref}}\in \mathcal{M}$, a time alignment procedure with respect to $M_{\textrm{ref}}$ is a map which to each motion $M$ in $\mathcal{M}$ associates a unique time warping $\varphi \in \operatorname{Diff}^{+}([0, 1])$ such that $M_{\textrm{ref}}$ and $M\circ\varphi$ are visually as closed as possible. 

\subsection{Properties of a temporal alignment procedure}\label{properties}
Intuitively, a temporal alignment procedure should satisfy:
\begin{property}\label{Property1a}
The optimal diffeomorphism aligning $M_{\textrm{ref}}$ with respect to itself should by the identity map $\textrm{id}: [0, 1]\rightarrow [0, 1]$, $t\mapsto t$ (reflexivity);
\end{property}
\vspace{-0.3cm}
\begin{property}\label{Property1b}
if $M_2$ is visually as close as possible to $M_1$ than $M_1$ should be visually as close as possible to $M_2$ (symmetry);
\end{property}
\vspace{-0.3cm}
\begin{property}\label{Property1c}
The optimal diffeomorphism aligning $M_{\textrm{ref}}\circ \varphi$ with respect to $M_{\textrm{ref}}$ should be $\varphi^{-1} \in \operatorname{Diff}^{+}([0, 1])$;
\end{property}
\vspace{-0.3cm}
\begin{property}\label{Property1d}
If two motions $M_1$ and $M_2$ are considered as visually as closed as possible, then for any $\varphi \in \operatorname{Diff}^{+}([0, 1])$, $M_1\circ\varphi$ and $M_2\circ \varphi$ should be visually as close as possible ($\operatorname{Diff}^{+}([0, 1])$-equivariance);
\end{property}
\vspace{-0.3cm}
\begin{property}\label{Property1e}
If $M_1$ is visually as close as possible to $M_2$ and $M_2$ is visually as close as possible to $M_3$, than $M_1$ should be visually as close as possible to $M_3$ (transitivity).
\end{property}

\subsection{Temporal Alignment procedures as projections on slices}\label{projections}

Properties~\ref{Property1a}, \ref{Property1b} and \ref{Property1e}  imply that ``being visually as close as possible'' is an equivalence relation $\sim$ on the set of motions $\mathcal{M}$. 
In this setting, given a reference motion $M_{\textrm{ref}}$, the set of motions visually as close as possible to $M_{\textrm{ref}}$ is called the equivalence class of $M_{\textrm{ref}}$.
Properties~\ref{Property1c} and ~\ref{Property1d} imply that the group of diffeomorphisms acts equivariently on the set of equivalence classes, i.e. for any $\varphi \in \operatorname{Diff}^{+}([0, 1])$,
$M_1 \sim M_2 \Leftrightarrow M_1\circ \varphi \sim M_2\circ \varphi.
$
The uniqueness of the optimal time warping of a motion with respect to a reference motion $M_{\textrm{ref}}$ implies that the orbit of any motion under the action of $\operatorname{Diff}^{+}([0, 1])$ by reparameterizations intersects the equivalence class of $M_{\textrm{ref}}$ at a unique point. Mathematically this means that the equivalence class of  a reference motion  is a global slice to the set of $\operatorname{Diff}^{+}([0, 1])$-orbits, and that a temporal alignment procedure provides a $\operatorname{Diff}^{+}([0, 1])$-invariant projection on it.

\vspace{-0.2cm}
\section{Experimental results}\label{results}
In this section, we first present some alignment procedures that will be used in the paper (Subsection~\ref{procedures}). Next, based on our mathematical formulation, we present a consistency test and compare the accuracy  of the alignment procedures under consideration (Subsection~\ref{Consistency_check_section}).
We notice that the coarse alignment procedure given by keyframe correspondences has an overall good performance for very low computational cost. For this reason, we incorporate keyframe correspondences into the dynamic programming algorithm  in order to improve the other alignment procedures (Subsection~\ref{incorporating}). This gain can be explained by the fact that the group of invariances of motions under consideration is $\mathbb{R}\times \operatorname{SO}(2)$ instead of $\operatorname{SE}(3)$, i.e. variations in the elevation of joints contain crucial information for synchronizing two motions. Reincorporating this information, even in a coarse manner, improves the performance of algorithms using $\operatorname{SE}(3)$-invariant features.
\vspace{-0.2cm}
\subsection{Examples of Alignment procedures}\label{procedures}
\subsubsection{Temporal alignment using keyframes.}\label{height_section}
\begin{figure*}[ht!]
\includegraphics[width=0.3\textwidth]{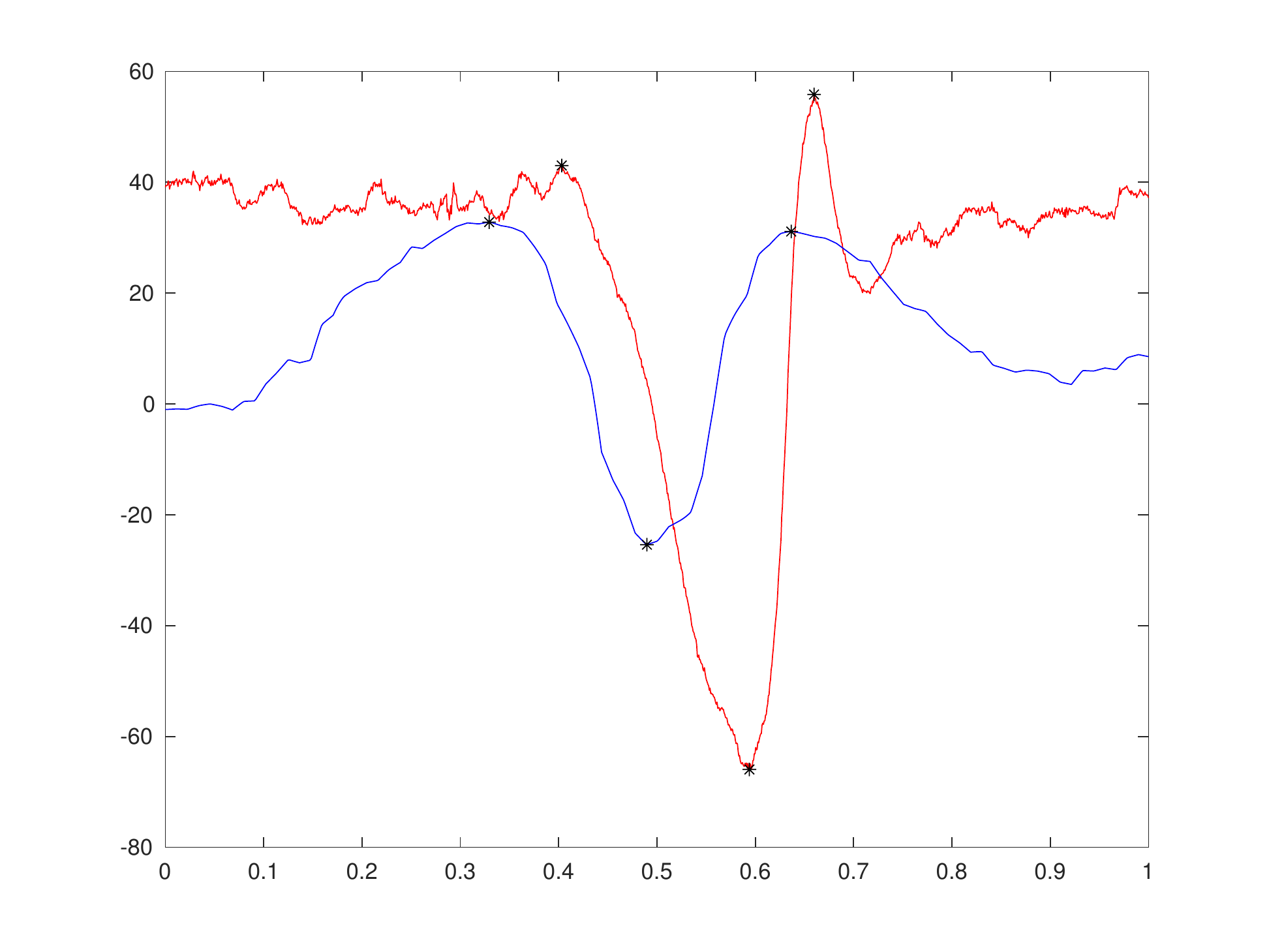}
\includegraphics[width=0.3\textwidth]{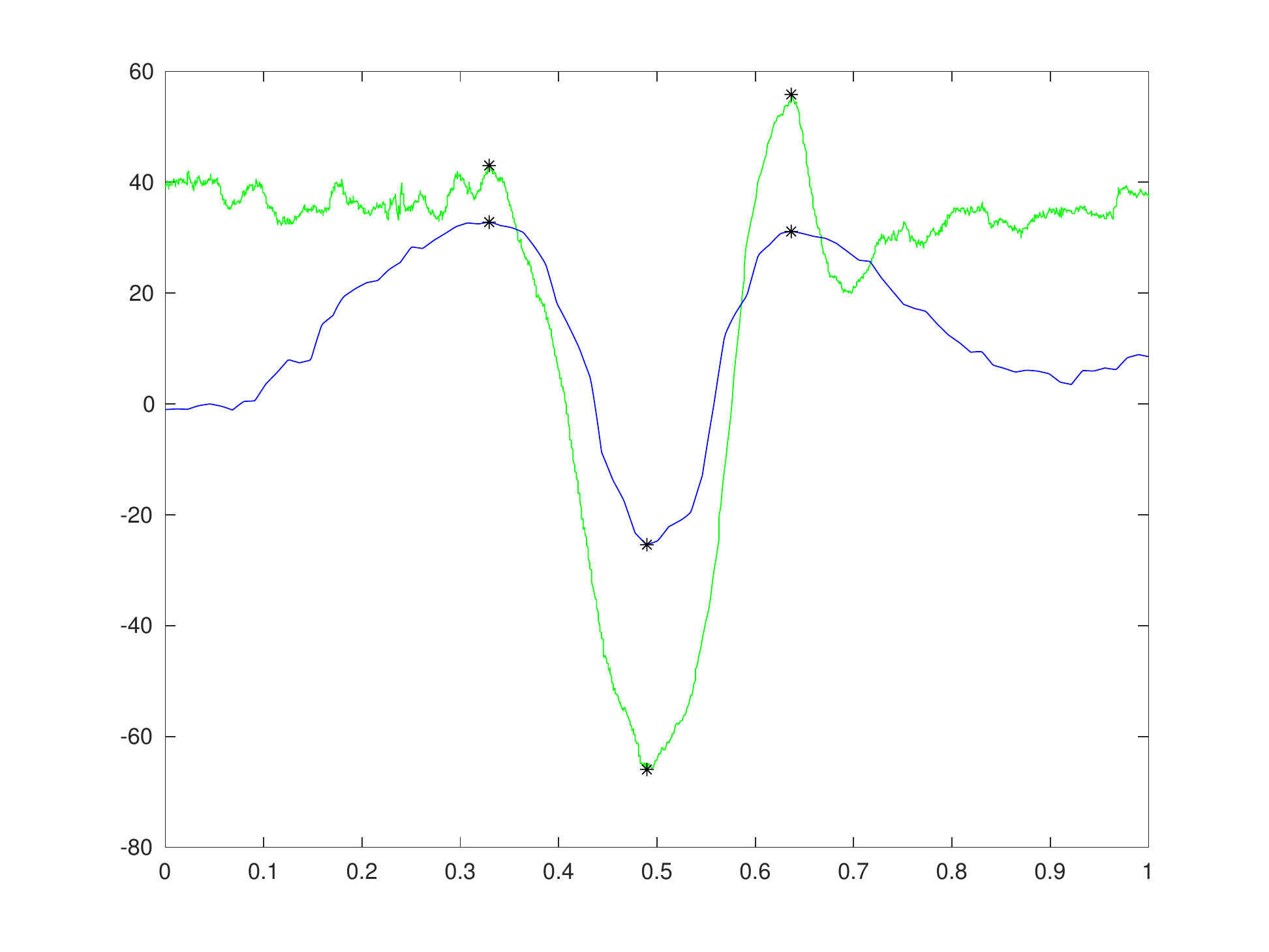}
\includegraphics[width=0.3\textwidth]{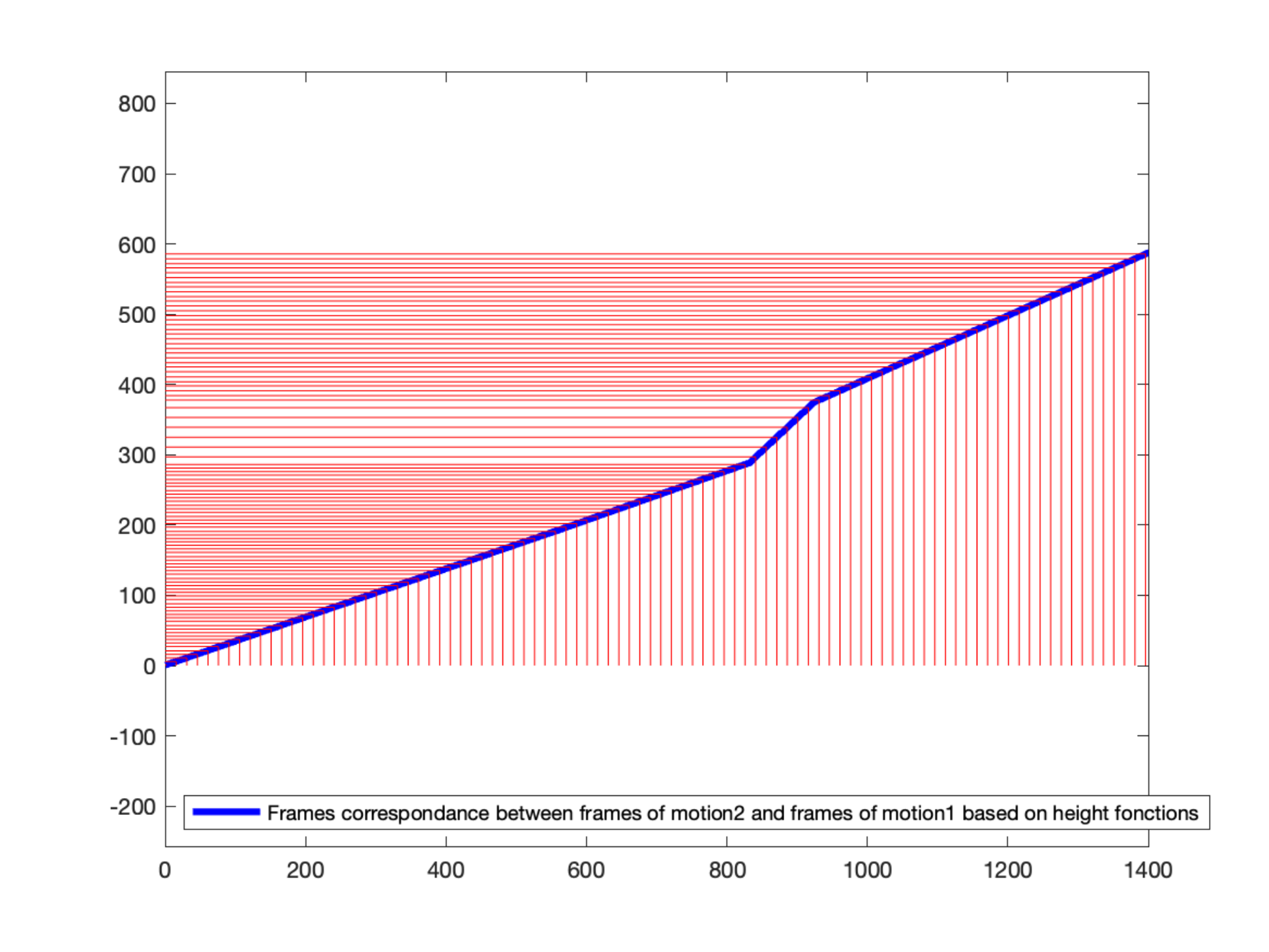}
\caption{\scriptsize{Elevation of the arm of a tennis player for two different motions before alignment (left) and after alignment (middle), and the frame correspondence aligning keyframes.}\label{keyframe_correspondence}}
\end{figure*}

As a coarse alignment procedure, we have implemented  a keyframe correspondence based on the elevation of the arm holding the racket. It allows to give a temporal bounding box around the movement of interest. 
For each selected joint of the arm holding the racket, the algorithm detects 3 keyframes:
\vspace{-0.2cm}
\hspace{1cm}\begin{enumerate}
\item the first frame with the highest $z$-coordinate of the joint;
\item the first frame with the lowest $z$-coordinate of the joint;
\item the second frame with the highest $z$-coordinate of the joint.
\end{enumerate}
\vspace{-0.2cm}
A frame correspondence between two motions is then calculated as the piecewise-linear frame correspondence mapping keyframes to keyframes (see Fig.~\ref{keyframe_correspondence}).

\vspace{-0.3cm}
\subsubsection{Temporal alignment using SRVT on trajectories of joints.}\label{SRVT_section}
We will use a transformation called Square Root Velocity Transform (SRVT) that allows temporal alignment of 3D-curves \cite{Srivastava}. 
For this alignment procedure, Dynamic Programming is used to minimize the $L^2-$distance between the SRVT-transforms of the trajectories of a selection of joints. 

\vspace{-0.3cm}
\subsubsection{Temporal alignment using Gram-matrices.}\label{Gram_section}
In this method, we compute  the Gram-matrices \cite{Celozzi,Kacem} associated to the joints positions of each skeleton and align them in the space of positive semi-definite matrices using Dynamic Programming.
In our experiments, we used 10 active joints, namely \textit{``Ankle left'', ``Ankle right'', `` Hip left'',  ``Hip right'', ``Knee left'', ``Knee right'', ``Spine low'', ``Spine high'', ``Racket hand'', ``Racket top''}.

\vspace{-0.3cm}
\subsubsection{Temporal alignment using curves on the group of rotations $\operatorname{SO}(3)$.}
We construct the moving frame associated to the trajectory of a joint, which is a curve on the group of rotations $\operatorname{SO}(3)$ and we reparameterize it in a canonical way.  
We use Dynamic Programming to align the curvature and torsion functions, which corresponds to a cost function measuring the difference of velocities of the curves on $\operatorname{SO}(3)$. 

\vspace{-0.3cm}
\subsubsection{Temporal alignment using curves on the sphere $\mathbb{S}^2$.}\label{section_S2}
Given a motion, we create a curve on the sphere $\mathbb{S}^2$ by joining a given joint of the skeleton to the center of the body, and by normalizing the vector obtained. Given two curves on the sphere corresponding to two different motions of a joint, we align them using Dynamic Programming where the cost function uses the SRV transform for homogeneous spaces introduced by Celledoni el al. \cite{Celledoni2,Celledoni3}.

\vspace{-0.2cm}
\subsection{Combining time wrapings for different joints}

Except the alignment procedure based on Gram-matrices, all the alignment procedures presented in Subsection~\ref{procedures} compute one diffeomorphism per joint. In order to combine the results, we use a weighted average or a median. These methods are more efficient from a computational point of view then the global method using Gram-matrices, because the joints calculations can be parallelized.

\vspace{-0.2cm}
\subsection{Consistency check}\label{Consistency_check_section}
\begin{figure*}[ht!]
\includegraphics[width=0.78\textwidth]{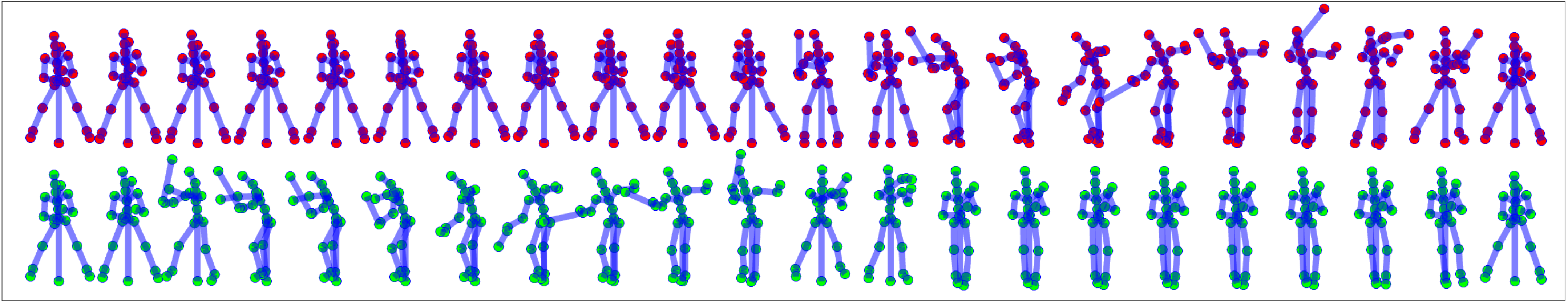}
\includegraphics[width=0.21\textwidth]{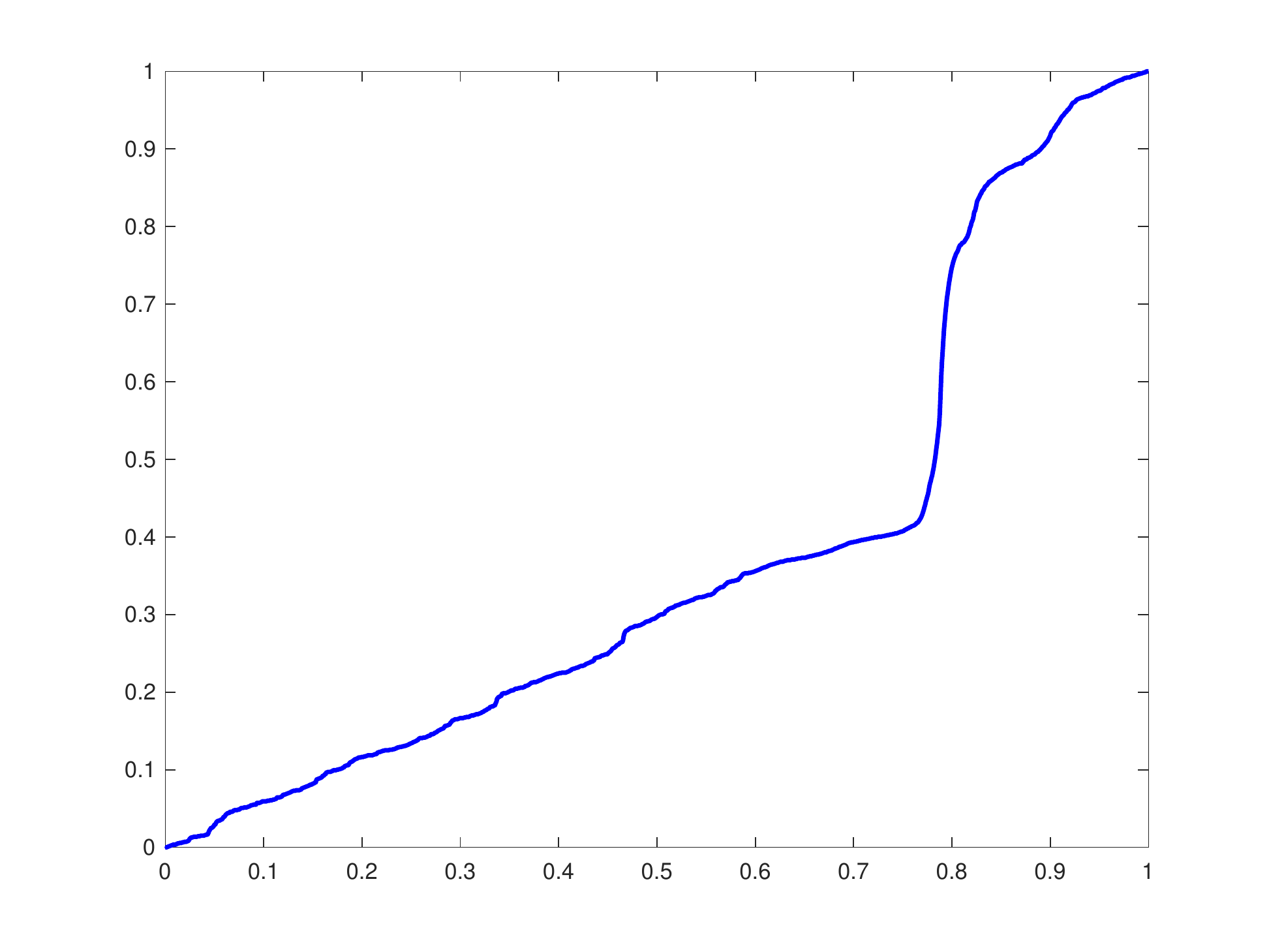}
\caption{\scriptsize{Bottom left: initial motion, upper line left: artificially reparameterized motion, Right: reparameterization applied.}\label{skeletons_reparameterization_figure}}
\end{figure*}

We designed a consistency check for testing the accuracy of each implemented algorithm. Namely, according to Property~\ref{Property1c}, each alignment procedure to a reference motion $M_{\textrm{ref}}$ taking $M_{\textrm{ref}}\circ\varphi$ as input should give as output $\varphi^{-1}$ for any reparameterization $\varphi\in \operatorname{Diff}^{+}([0, 1])$.
To test this, we have implemented a function which takes a skeleton motion and an arbitrary reparameterization as input, and creates a new skeleton motion given by the frame correspondence applied to the initial skeleton. The result is a skeleton motion that performs exactly the same action but with a different rate. 
In Fig.~\ref{skeletons_reparameterization_figure},  
one can see an example of original motion at the bottom line (joints in green), and the same motion artificially reparameterized at the top line (joints in red) where we can clearly see that the movement starts later than in the initial motion. The reparameterization applied is displayed on the right of the same Figure.
\vspace{-0.6cm}
\vspace{-0.6cm}

\begin{figure*}[ht!]
\includegraphics[width = 0.5 \linewidth]{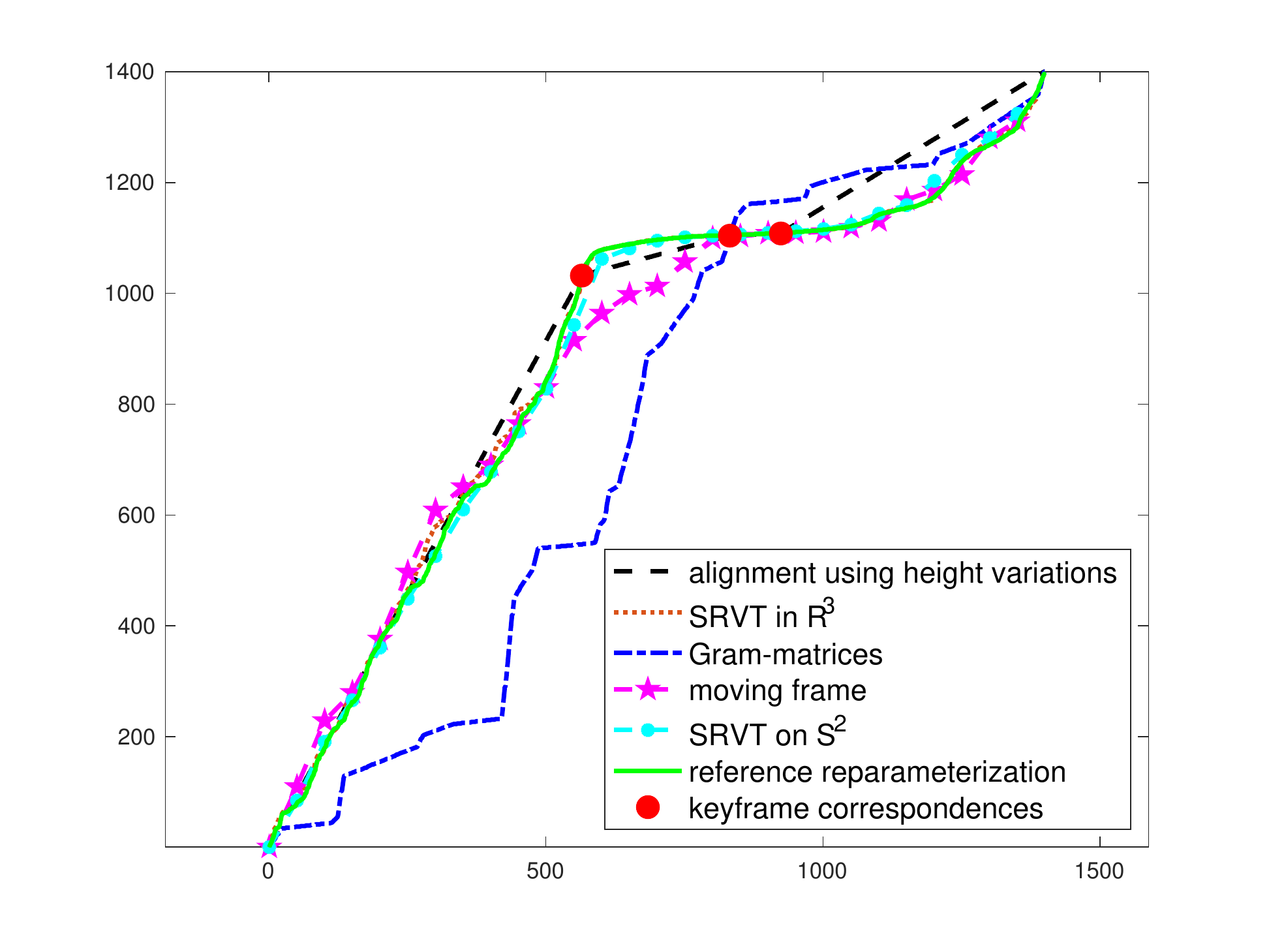} 
\includegraphics[width = 0.5 \linewidth]{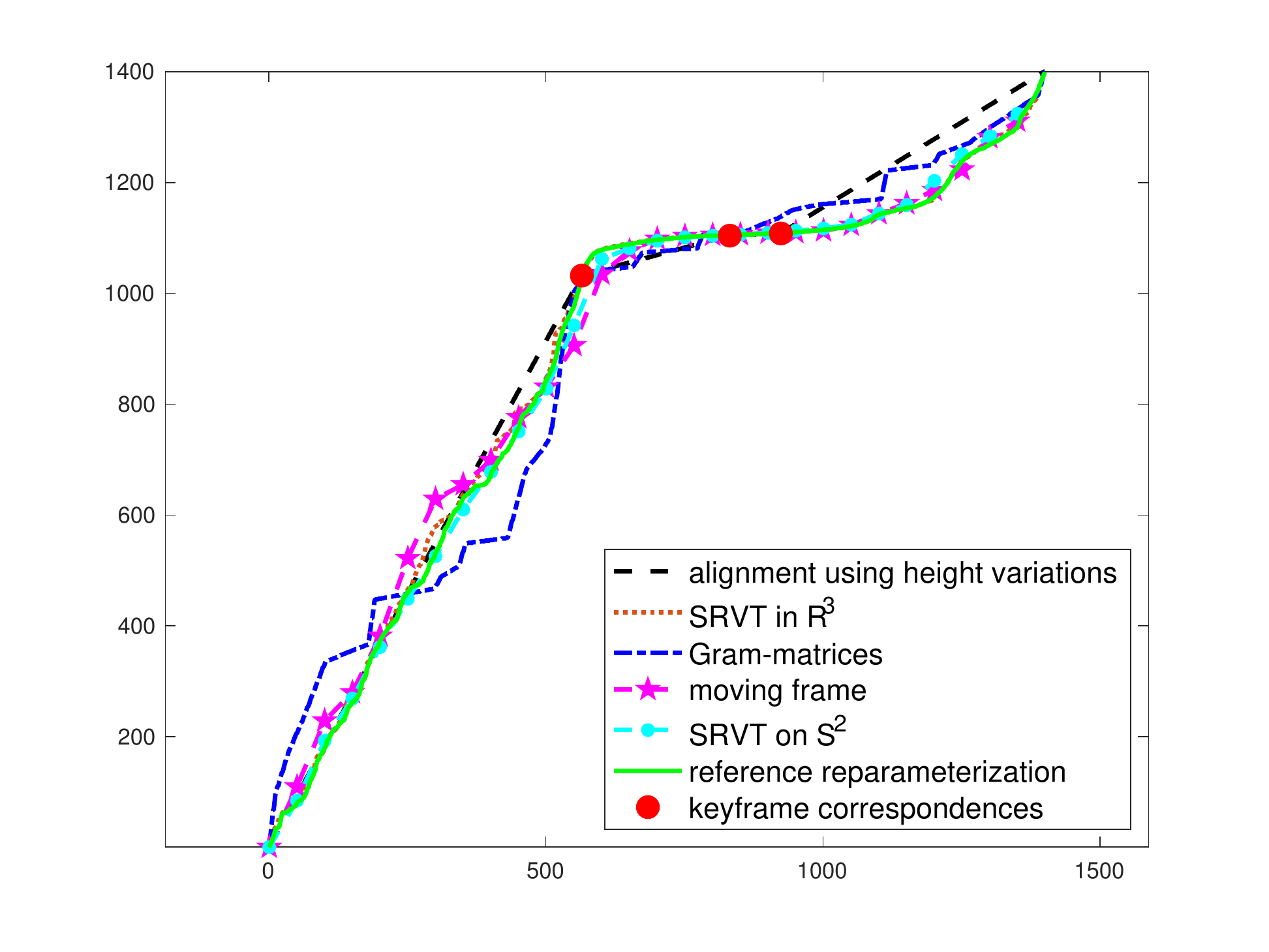} 
\caption{\scriptsize{Frame correspondence obtained with different alignment procedures without keyframe correspondence (left) and with keyframe correspondence (right). Coarse alignment based on elevation of arm (black dashed line), SRVT on $\mathbb{R}^3$ (orange dotted line), Gram-matrices (blue Dash-dotted line), moving frames (magenta dashed line with stars), SRVT on $\mathbb{S}^2$ (cyan dashed line with discs), and reference correspondence (green solid line). }\label{Consistency_figure}}
\end{figure*}

We have tested the consistency of each algorithm presented in Section~\ref{procedures}. An example of outputs corresponding to the alignment of the upper motion (red joints) of Fig.~\ref{skeletons_reparameterization_figure} to the initial motion (green joints)
as well as the ground-truth provided by the inverse of the diffeomorphism given in Fig.~\ref{skeletons_reparameterization_figure} can be seen in Fig.~\ref{Consistency_figure} left. The right picture in Fig.~\ref{Consistency_figure} corresponds to an improvement of the algorithms explained in next section. The mean $L^1$-errors between the output of each alignment procedure and the ground-truth computed over 7 experiments with varying number of frames is recorded in Table~\ref{L1}, as well as average computational times under the same conditions.
\vspace{-0.6cm}
\vspace{-0.6cm}
\begin{table}[ht!]
\centering
\caption{\scriptsize{mean $L^1$-errors between the frame correspondence provided by each alignment procedure and the ground truth over 7 experiments with number of frames between 50 and 185, and average computational times on a Macbook M1. Left: using dynamic programming (DP), right: using anchored dynamic programming (ADP).} \label{L1}}
\begin{tabular}{|l|r|r|r|r|}
\hline  & Error with DP & Time with DP & Error with ADP & Time with ADP \\
\hline
$L^1$-Error SRVT in $\mathbb{R}^3$ &   \textbf{0.95\%} & 23s 227ms   &  \textbf{0.93\%} & 7s 692ms \\
$L^1$-Error Gram-matrices &   16.30\% &9m 38s 706ms  &   10.20\% &2m 49s 257ms \\
$L^1$-Error Moving Frames & 6.32\% & 23s 771ms  & 2.64\% & 8s 184ms\\
$L^1$-Error SRVT  on $\mathbb{S}^2$ &  1.15\% & 23s 250ms &   1.12\% & 7s 640ms\\
\hline
\end{tabular}
Baseline: $L^1$-Error alignement of keyframes = $4.57\%$, computational time = \textbf{5ms}.
\end{table}

\vspace{-0.5cm}
\subsection{Incorporating keyframe correspondences into Dynamic Programming}\label{incorporating}
\begin{figure*}[ht!]
\centering

\includegraphics[width = 0.2 \linewidth, angle = 90]{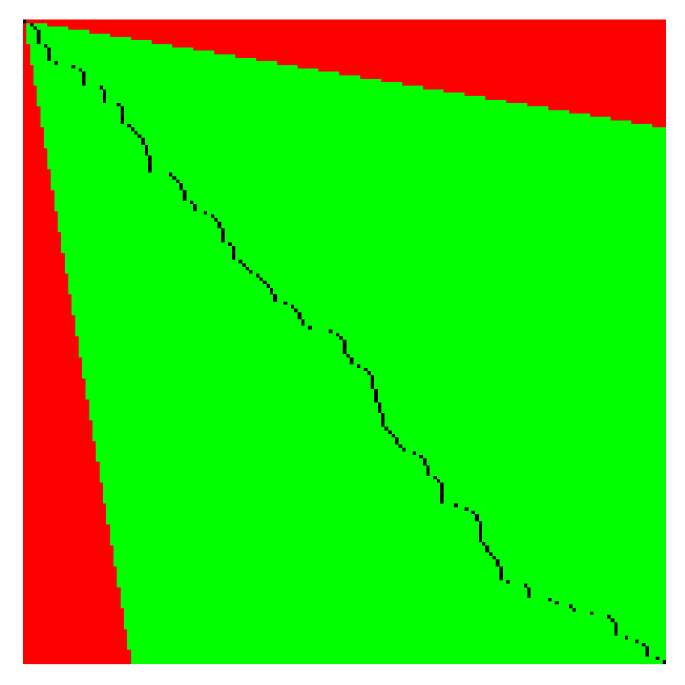} 
\includegraphics[width = 0.2 \linewidth, angle = 90]{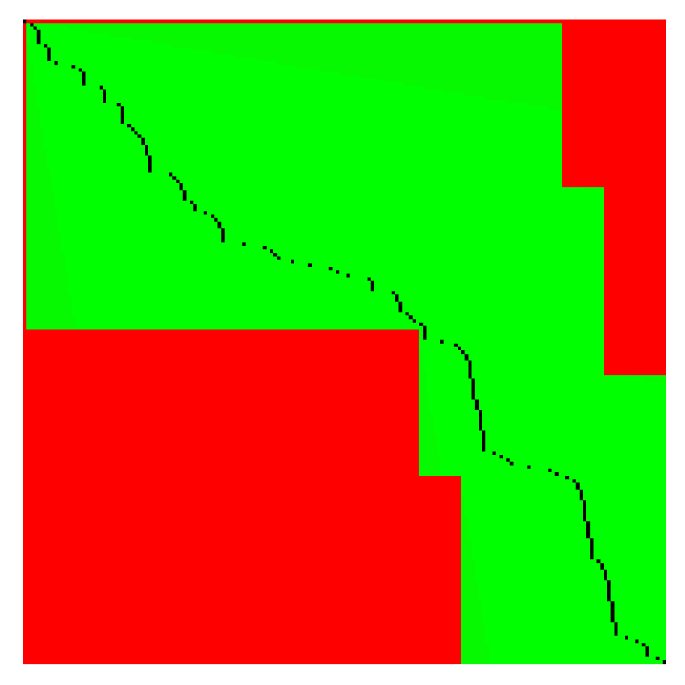} 
\includegraphics[width = 0.2 \linewidth, angle = 90]{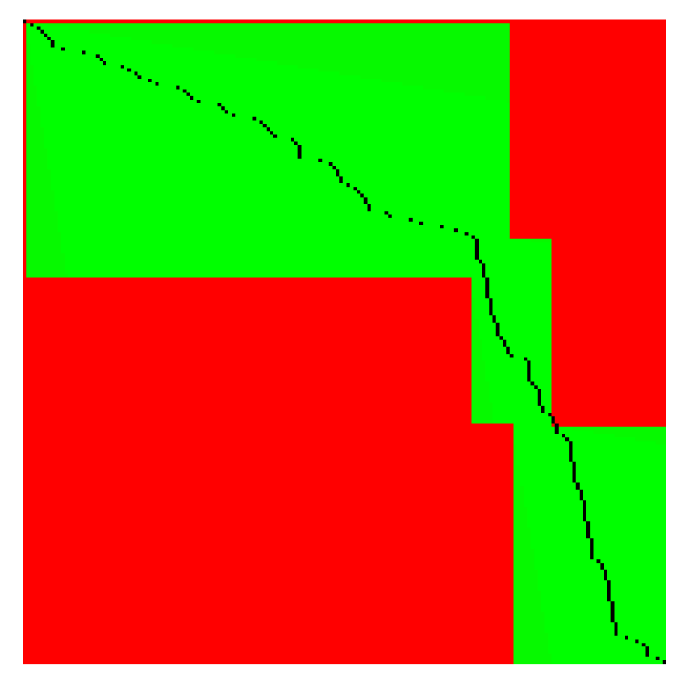} 
\includegraphics[width = 0.2 \linewidth, angle = 90]{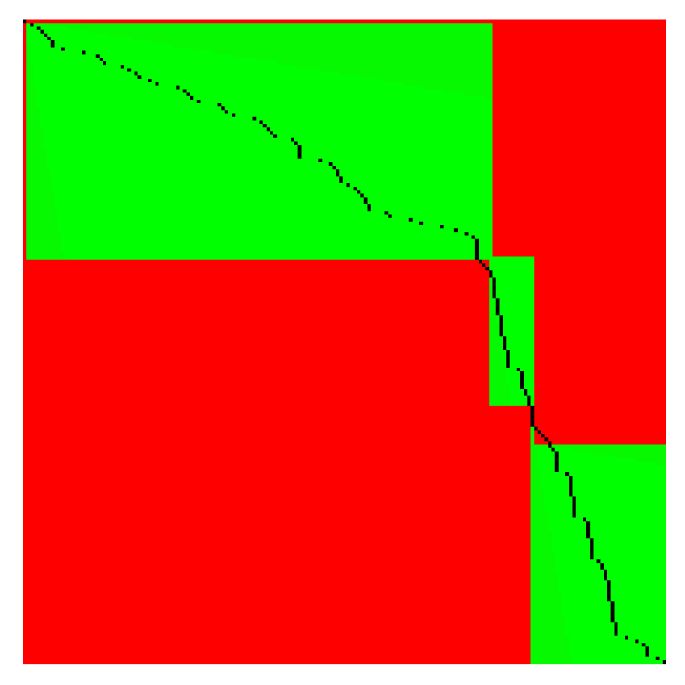} 
\caption{\scriptsize{Energy profile for dynamic programming and example of path of lowest energy. Left: classical dynamic programming, second and third from left: anchored dynamic programming around the nodes provided by the coarse alignment procedure based on keyframes with tolerance of 1/4  and 1/20 of total number of frames respectively, right: anchored dynamic programming with minimal tolerance.}\label{anchoredDP}}
\end{figure*}

\begin{figure*}[ht!]
\centering
\includegraphics[width=\textwidth]{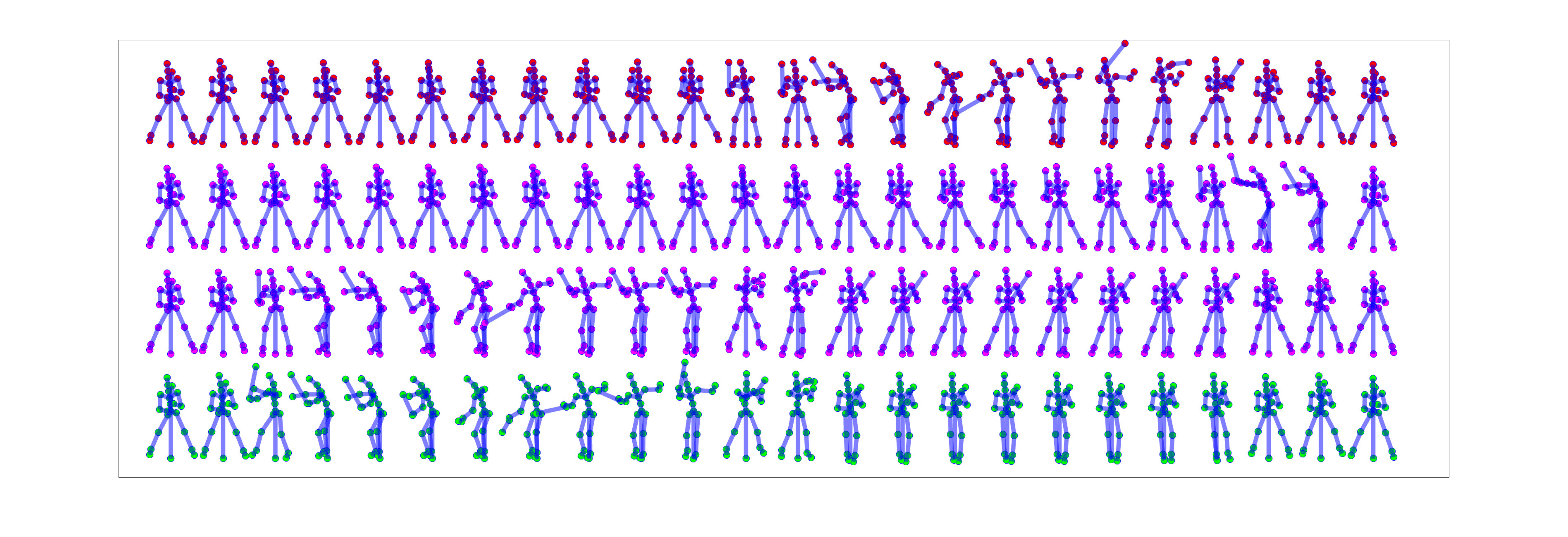}
\caption{\scriptsize{Bottom: initial motion (ground truth), upper line: artificially reparameterized motion obtained in Section~\ref{Consistency_check_section} given as input, second line: aligned motion using moving frames on 100 frames with classical dynamic programming, third line: aligned motion using moving frames on 100 frames with anchored dynamic programming with zero tolerance.}\label{mv_anchored}}
\end{figure*}

In order to take benefit of the stable good performance with low computational cost provided by the coarse alignment procedure based on the elevation of the arm holding the racket, we have modified the dynamic programming algorithm  to incorporate keyframe correspondences. Each keyframe correspondence can be thought as a node that should be traversed by the optimal time warping.
The resulting anchored dynamic programming finds the path of minimal energy in a landscape that is shaped according to a desired tolerance around each node as in Fig.~\ref{anchoredDP}. The modified alignment procedures are computationally more efficient (less nodes to visit) and more accurate (see Table~\ref{L1}). An example of aligned motions by the procedure using moving frames with classical dynamic programming and with anchored dynamic programming is displayed in Fig.~\ref{mv_anchored}.
\vspace{-0.6cm}
\vspace{-0.6cm}
\vspace{-0.6cm}

\section{Conclusion}

We gave a mathematical formulation of the task consisting of synchronizing human motion data from multiple recordings, which allowed us to elaborate a test to check the consistency of any temporal alignment procedure. In order to include the information gained by a coarse alignment procedure based on keyframes in any method, we implemented a variant of dynamic programming ensuring that associated keyframes are in correspondence. For this algorithm, each keyframe correspondence creates a node by which the optimal time warping has to pass, like a boat that needs to drop anchor in a port. The improvement of the temporal alignment procedures by the anchored dynamic programming can be explained by the fact that the group of symmetries of extracted features mismatched the group of symmetries of the motions under consideration. The lost information was partially recovered by forcing keyframe correspondence. At the same time the complexity of the algorithms decreased significantly.

\subsubsection{Acknowledgements} We thank VR Motion Learning GmbH \& Co KG for providing us with their dataset of tennis motions. The first author is supported by FWF grant I 5015-N. This work has been funded by grant F77 of the Austrian Science Fund FWF (SFB "Advanced Computational Design", SP5).
%
%
%
%

\end{document}